\documentclass[11
pt,a4paper]{article}
\usepackage[utf8]{inputenc}
\usepackage[english]{babel}
\usepackage[T1]{fontenc}
\usepackage{amsmath}
\usepackage{amsfonts}
\usepackage{amssymb}
\usepackage{amsthm} 
\usepackage[all]{xy} 
\usepackage[a4paper]{geometry}
\usepackage{mathrsfs}
\usepackage{frcursive}
\usepackage{color}
\usepackage{tikz-cd}
\usepackage{mathtools}
\usepackage{graphicx}

 \theoremstyle{plain}
 \newtheorem{thm}{Theorem}[section]
 
 \newtheorem{prop}[thm]{Proposition}
 \newtheorem{lem}[thm]{Lemma}
 \newtheorem{coro}[thm]{Corollary}

 \theoremstyle{remark}
 \newtheorem{rem}[thm]{Remark}

\author{Sylvain Gaulhiac}
\title{Comparison between admissible and de Jong coverings of rigid analytic spaces in mixed characteristic}

\date{}

 \font\bf= cmbx10 at 10pt

 \newcommand{\A}{\mathbb A}
 \newcommand{\C}{\mathbb{C}}
 
 \newcommand{\G}{\mathbb{G}}

 \newcommand{\N}{\mathbb{N}}
 \newcommand{\PP}{\mathbb{P}}
 \renewcommand{\P}{\PP}
 \newcommand{\Q}{\mathbb{Q}}
 \newcommand{\R}{\mathbb{R}}

 \newcommand{\Aut}{\mathrm{Aut}}

\newcommand{\an}{\mathrm{an}}

 \newcommand{\CC}{\mathcal{C}}

 \newcommand{\iso}{\xrightarrow{\sim}}

\begin{document}

\maketitle

\section*{Introduction}

In the two recent works \cite{ALY1} and \cite{ALY2}, Achinger, Lara and Youcis made substantial progress in the theory of covering spaces in non-archimedean analytic geometry. If $k$ is a complete non-archimedean field, and $X$ an adic space locally of finite type over $\mathrm{Spa}(k)$, they define a new class of coverings, the \emph{geometric coverings}. This class is much larger than the class of tempered coverings defined by André in \cite{And} and even larger than the one studied by de Jong in \cite{DJg}. It is closed under disjoint unions, étale local on $X$ and forms a tame infinite Galois category when $X$ is connected. Let's explain shortly the purpose of this note :

\bigskip

We introduce some notations (from \cite{ALY1}).  If $X$ is an adic or Berkovich space, $\textbf{Ét}_X$ (resp. $\textbf{FÉt}_X$) is the category of étale (resp. finite étale) objects over $X$, and $\textbf{UFÉt}_X$ stands for the category of disjoint unions of objects of $\textbf{FÉt}_X$. By $\textbf{Cov}_{X}^{\mathrm{oc}}$ we mean the full subcategory of $\textbf{Ét}_X$ consisting of étale morphisms that are locally for the overconvergent topology (i.e. for the Berkovich topology of the corresponding Berkovich space) disjoint union of finite étale coverings. These morphisms will be called \emph{de Jong coverings}, since they were first defined and studied by de Jong in his article \cite{DJg}. Similarly, $\textbf{Cov}_{X}^{\mathrm{adm}}$ stands for the full subcategory of $\textbf{Ét}_X$ of morphisms that are locally for the admissible topology of $X$ disjoint union of finite étale coverings. These categories are full subcategories of the category $\textbf{Cov}_X$ of \emph{geometric coverings} of $X$ defined in \cite{ALY1}. They are not closed under disjoint unions, but this can be fixed by defining the category $\textbf{UCov}_{X}^{\tau}$ of disjoint unions of objects of $\textbf{Cov}_{X}^{\tau}$, for $\tau\in \lbrace \mathrm{oc}, \mathrm{adm}\rbrace$. Achinger, Lara and Youcis showed ($5.2.1$, \cite{ALY1}) that these categories $\textbf{UCov}_{X}^{\tau}$ have good stability properties : when $X$ is connected, any geometric point $\overline{x}$ induces a natural fiber functor $F_{\overline{x}}$ such that  the pair $(\textbf{UCov}_{X}^{\tau}, F_{\overline{x}})$ defines a \emph{tame infinite Galois category} in the sense of \cite{BS15}. In particular, when $X$ is connected, $(\textbf{UCov}_{X}^{\tau}, F_{\overline{x}})$ has a fundamental group $\pi_1^{\mathrm{dJ}, \tau}(X,\overline{x})$ called the \emph{$\tau$-adapted de Jong fundamental group} of the pair $(X,\overline{x})$. It is a Noohi group ($7.1$, \cite{BS15}), and the functor $$F_{\overline{x}} : \textbf{UCov}_{X}^{\tau}\to \pi_1^{\mathrm{dJ}, \tau}(X,\overline{x})-\textbf{Set}$$ is an equivalence of categories. \\
 
 As the admissible topology on $X$ is finer than the overconvergent Berkovich topology, there is a natural inclusion $\textbf{Cov}_{X}^{\mathrm{oc}}\subseteq \textbf{Cov}_{X}^{\mathrm{adm}}$, inducing a morphism of Noohi groups : $\pi_1^{\mathrm{dJ, \, adm}}(X)\to \pi_1^{\mathrm{dJ, \,oc}}(X)$ with dense image. A question initially asked by de Jong (\cite{DJg}, p. $106$) is whether the inclusion $\textbf{Cov}_{X}^{\mathrm{oc}}\subseteq \textbf{Cov}_{X}^{\mathrm{adm}}$ is strict, in other terms whether the property of being a de Jong covering can be checked locally for the admissible topology. Achinger, Lara and Youcis partially answered this question, showing that : 
\begin{itemize}
\item[•]the inclusion $\textbf{Cov}_{X}^{\mathrm{oc}}\subseteq \textbf{Cov}_{X}^{\mathrm{adm}}$ is strict when $k$ is of finite characteristic (\cite{ALY1}, $5.3.4$),
\item[•]$\textbf{Cov}_{X}^{\mathrm{oc}}= \textbf{Cov}_{X}^{\mathrm{adm}}$ when $k$ is discretely valued of equal characteristic $0$, and $X$ smooth, quasi-paracompact and quasi-separated (\cite{ALY2}, $4.17$).
\end{itemize}

In this note, I show that the inclusion $\textbf{Cov}_{X}^{\mathrm{oc}}\subseteq \textbf{Cov}_{X}^{\mathrm{adm}}$ is strict in the mixed characteristic case $(0,p)$, when $k$ is $p$-closed (i.e. contains all the $p^{\mathrm{th}}$-roots of all the elements). The arguments are very similar to the ones used in (\cite{ALY1}, $5.3.4$). The idea is to exhibit an admissible covering of an annulus $\CC$ that is not a de Jong covering, built by taking the disjoint union and gluing in a clever way of some sequences of $\mu_p$-torsors on two intersecting sub-annuli whose splitting domain is shrinking such that the overconvergent topology is not fine enough to control the limit of this shrinking. The contribution of this note is to specify some good $\mu_p$-torsors adapted to the mixed characteristic case that are replacing the Artin-Schreier coverings used in the finite characteristic case.

\bigskip
As a corollary, we get that if $\overline{x}$ is a geometric point of $\CC$, the morphism of fundamental groups $$\pi_1^{\mathrm{dJ, \, adm}}(\CC, \overline{x})\to \pi_1^{\mathrm{dJ, \,oc}}(\CC,\overline{x}) $$ is not an isomorphism.

\bigskip
\textbf{Remark :} In \cite{ALY1}, the authors  work with an adic space $X$ locally of finite type over $\mathrm{Spa}(k)$ instead of the corresponding Berkovich space $X^{\mathrm{Berk}}$, which can be characterised topologically as the universal separated quotient of $X$ (see \cite{ALY1}, 1.2.2, or \cite{FK}). This framework of adic spaces enable them, when defining geometric coverings, to think about partial properness in term of a valuative criterion instead of the notion of boundaryless in Berkovich langage. However, in this note, without the necessity of working with the notion of partial properness, I will rather work directly with Berkovich spaces, where the notion of the skeleton of a curve as well as the canonical retraction on the skeleton find its natural framework. The terminology of Berkovich spaces used here (points of the form $\eta_{a,r}$, type-$2$ points, skeleton $S^\an(X)$) is classical, see for instance (\cite{Gau}, $\S 2$).

\section{Splitting conditions of some $\mu_p$-torsors}

Let $k$ be a $p$-closed (i.e. contains all the $p^{\mathrm{th}}$-roots of all the elements) complete non-archimedean field of mixed characteristic $(0,p)$ : $\mathrm{char}(k)=0$ and $\mathrm{char}(\widetilde{k})=p$, where $\widetilde{k}$ is the residue field of $k$. We assume that the absolute value on $k$ is normalized such that $\vert p\vert=p^{-1}$. The field $\C_p=\widehat{\overline{\Q_p}}$ with the usual $p$-adic absolute value is an example of such field. 

\begin{lem}\label{lemme sur la distance entre deux racines de l'unité} 
Let $\xi$ and $\xi'$ be two distinct $p^{\mathrm{th}}$-roots of unity in $k$. Then $\vert \xi - \xi'\vert= p^{-\frac{1}{p-1}}$.
\end{lem}

\begin{proof}
Write $\displaystyle{\Phi_p=\frac{X^p-1}{X-1}=\sum_{i=0}^{p-1}X^i=\prod_{\xi\in \mu'_p}X-\xi\in \Q[X]}$ for the $p^{\mathrm{th}}$ cyclotomic polynomial, where $\mu_p'$ stands for the set of the $p-1$ primitive $p^{\mathrm{th}}$-roots of unity in $k$. The evaluation at $1$ gives : $p=\prod_{\xi\in \mu'_p}1-\xi.$ For $\xi$ describing $\mu'_p$, all the $1-\xi$ have the same norm since they are on the same $\Aut(k/\Q_p)$-conjugacy class. Therefore, $\vert 1-\xi \vert=p^{-\frac{1}{p-1}}$, and we obtain the result since multiplication by any $p^{\mathrm{th}}$-root of unity is an isometry of $k$. 
\end{proof}

An étale morphism $\varphi : Y\to X$ between two $k$-analytic Berkovich curves \emph{totally splits} over a point $x\in X$ if for any $y\in \varphi^{-1}(\lbrace x\rbrace)$, the extension of completed residual fields $\mathscr{H}(x)\to \mathscr{H}(y)$ is an isomorphism. When $\varphi$ is of degree $n$, $\varphi$ totally splits over $x$ if and if the fibre $\varphi^{-1}(\lbrace x\rbrace)$ has exactly $n$ elements, which is the same as saying that locally, over a neighbourhood of $x$, $\varphi$ is a topological covering (it follows easily from \cite{Ber}, Lemma $3.4$).
The following proposition describes the splitting sets of the $\mu_{p^h}$-torsor given by the function $\sqrt[p^h]{1+T}$.

\begin{prop}\label{décomposition du torseur sauvage} 
If $h\in \N^\times$, the étale covering $f=\G_m^\an \xrightarrow{z\mapsto z^{p^h}} \G_m^\an$ totally splits over a point $\eta_{z_0, r}$ satisfying $r\leqslant\vert z_0\vert=:\alpha$ if and only if : $r<\alpha p^{-h-\frac{p}{p-1}}$. More precisely, the inverse image of $\eta_{z_0, r}$ contains :
\begin{itemize}
\item[•]only one element when $r\in [\alpha p^{-\frac{p}{p-1}}, \alpha];$
\item[•]$p^i$ elements when $r\in[\alpha p^{-i-\frac{p}{p-1}}, \alpha p^{-i-\frac{1}{p-1}}[,$ with $1\leqslant i\leqslant h-1$;
\item[•]$p^h$ elements when $r\in[0, \alpha p^{-h-\frac{1}{p-1}}[.$
\end{itemize}
\end{prop}

\begin{proof}
Let $z_1\in k^\times$ and $\rho\in \R_{\geqslant 0}$ satisfying $\rho\leqslant\vert z_1 \vert$.\
In order to compute $f(\eta_{z_1,\rho})$, notice that for any polynomial $P\in k[T]$:

\begin{align*}
\vert P\left(f(\eta_{z_1, \rho})\right) \vert &= \vert (P\circ f)(\eta_{z_1, \rho})\vert= \vert P(T^p)(\eta_{z_1, \rho})\vert\\
&=\sup_{x\in \mathrm{B}(z_1, \rho)}\vert P\circ f(x)\vert\\
&=\sup_{y\in f(\mathrm{B}(z_1, \rho))} \vert P(y)\vert
\end{align*}

As $k$ is $p$-closed, there exists $\widehat{\rho}>0$ such that $f(\mathrm{B}(z_1, \rho))=\mathrm{B(z_1^p,\widehat{\rho}\,)}$, which gives $f(\eta_{z_1, \rho})=\eta_{z_1^p,\widehat{\rho}}$.\\

In order to compute $\widehat{\rho}$, notice that $\widehat{\rho}=\vert (T-z_1^p)(f(\eta_{z_1, \rho}))\vert=\vert (T^p-z_1^p)(\eta_{z_1, \rho})\vert$, and : 
$$T^p-z_1^p=\left((T-z_1)+z_1\right)^p-z_1^p=\sum_{i=1}^p \binom{p}{i}z_1^{p-i}(T-z_1)^i=\sum_{i=1}^p \gamma_i (T-z_1)^i ,$$
where $\gamma_i= \binom{p}{i}z_1^{p-i}$, with : 
$$
\vert \gamma_i \vert = \left\{
    \begin{array}{ll}
        1& \hbox{if } i=p \\
        p^{-1}\vert z_1 \vert ^{p-i} & \hbox{if } 1\leqslant i\leqslant p-1
    \end{array}
\right.
$$
Therefore, $\widehat{\rho}=\vert (T-z_1^p)(f(\eta_{z_1, \rho}))\vert=\max_{1\leqslant i \leqslant p} \lbrace\vert \gamma_i \vert \rho^i \rbrace=\max \lbrace \rho^p, \left( p^{-1}\rho^i \vert z_1 \vert^{p-i}\right)_{1 \leqslant i \leqslant p-1}  \rbrace $. Since we assumed $\rho\leqslant \vert z_1 \vert$, we get $\widehat{\rho}= \max \lbrace \rho^p, p^{-1}\rho\vert z_1 \vert^{p-1}\rbrace $, that is to say : 
$$
\widehat{\rho} = \left\{
    \begin{array}{ll}
         p^{-1}\rho \vert z_1 \vert^{p-1}& \hbox{if } \rho \leqslant \vert z_1 \vert p^{-\frac{1}{p-1}} \\
        \rho^p & \hbox{if } \rho \geqslant \vert z_1 \vert p^{-\frac{1}{p-1}}
    \end{array}
\right.
$$

Consequently : 

$$
f(\eta_{z_1, \rho})= \left\{
    \begin{array}{ll}
         \eta_{z_1^p,p^{-1}\rho \vert z_1 \vert^{p-1}}& \hbox{if } \rho \leqslant \vert z_1 \vert p^{-\frac{1}{p-1}} \\
        \eta_{z_1^p, \rho^p} & \hbox{if } \rho \geqslant \vert z_1 \vert p^{-\frac{1}{p-1}}
    \end{array}
\right.
$$

\begin{figure}
\includegraphics[width=1\textwidth]{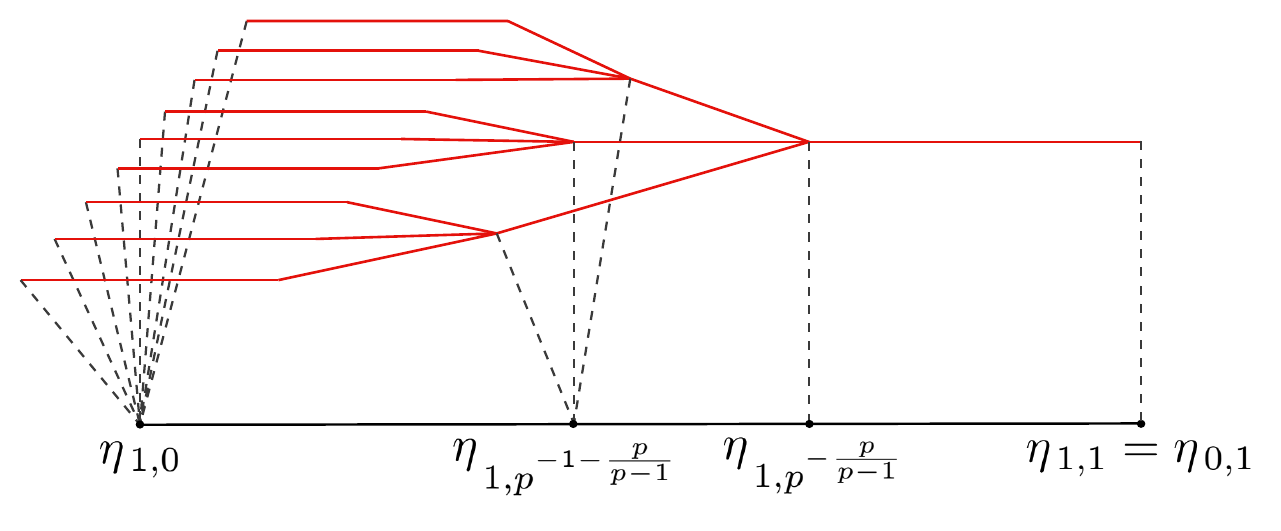}
\caption{Covering $\G_m^\an \xrightarrow{z\mapsto z^{9}} \G_m^\an$ with $p=3$, $h=2$ and $z_0=1$.}
\end{figure}

\bigskip
\item Let's try to find the preimages by $f$ of $\eta_{z_0, r}$, where $0\leqslant r\leqslant\alpha:=\vert z_0 \vert$. Define : 

$$
\widetilde{r} = \left\{
    \begin{array}{ll}
         rp \alpha^{-\frac{p-1}{p}}& \hbox{if } r \leqslant \alpha p^{-\frac{p}{p-1}} \\
        r^{\frac{1}{p}} & \hbox{if } r \geqslant \alpha p^{-\frac{p}{p-1}} 
    \end{array}
\right.
$$
From what is above, if $\widetilde{z_0}$ is a $p^{\mathrm{th}}$-root of $z_0$, then : $$\eta_{\widetilde{z_0}, \widetilde{r}}\in f^{-1}\left(\lbrace \eta_{z_0,r}\rbrace \right),$$ and $f^{-1}\left(\lbrace \eta_{z_0,r}\rbrace \right)$ consists of all conjugates $\eta_{\xi\widetilde{z_0}, \widetilde{r}}$ of $\eta_{\widetilde{z_0}, \widetilde{r}}$ for $\xi\in \mu_p$. Therefore : 
$$
f^{-1}\left(\lbrace \eta_{z_0,r}\rbrace \right) = \left\{
    \begin{array}{ll}
         \lbrace \eta_{\xi \widetilde{z_0},rp\alpha^{-\frac{p-1}{p}}}\rbrace_{\xi\in \mu_p} & \hbox{if } r \leqslant \alpha p^{-\frac{p}{p-1}} \\
         \lbrace \eta_{\xi \widetilde{z_0},r^{\frac{1}{p}}}\rbrace_{\xi\in \mu_p} & \hbox{if } r \geqslant \alpha p^{-\frac{p}{p-1}} 
    \end{array}
\right.
$$
Since $\vert \widetilde{z_0} \vert= \alpha^{\frac{1}{p}}$, we have $\vert \xi \widetilde{z_0}- \xi' \widetilde{z_0}\vert=\alpha^{\frac{1}{p}} p^{-\frac{1}{p-1}}$ as soon as $\xi\neq \xi' \in \mu_p$, from lemma \ref{lemme sur la distance entre deux racines de l'unité}. Thus, $f^{-1}\left(\lbrace \eta_{z_0,r}\rbrace \right)$ has a unique element if $r \geqslant \alpha p^{-\frac{p}{p-1}}$, $p$ otherwise. 
\item For the general case, with $h\geqslant 1$, a recursive reasoning on $h$ leads to the conclusion.

\end{proof}

\begin{rem}
A point $\eta_{z_0,r}\in \G_m^\an$ does not belong to the skeleton $S^\an(\G^\an_m)$ if and only if $r<\vert z_0 \vert$. Moreover, any type-$2$ point of the skeleton is of the form $\eta_{z_0,\vert z_0\vert}(=\eta_{0,\vert z_0\vert})$ for some $z_0\in \vert k^\times\vert$. Thus, from the density of type-$2$ points in the skeleton $S^\an(\G_m^\an)$, Proposition \ref{décomposition du torseur sauvage} implies in particular that the fiber of $f$ over any point of $S^\an(\G_m^\an)$ is a singleton.
\end{rem}\

If $\CC$ is a $k$-analytic annulus, let $\pi_\CC : \CC\twoheadrightarrow S^\an(\CC)$ be the canonical retraction on its skeleton.

\begin{coro}\label{corollary about splitting}

Let $\CC$ be a $k$-analytic annulus, considered as an analytic domain of $\A_k^{1,\an}$. Its skeleton $S^\an(\CC)$ is an interval (possibly finite, infinite, open, compact,\ldots) with extremities $x, y\in\P_k^{1,\an}$. Let $z\in ]x,y[$ a type-$2$ point of the interior of $S^\an(\CC)$. Then there exists a $\mu_p$-analytic torsor $f\in H^1(\CC,\mu_p)$ such that : 
\begin{itemize}
\item[•]$f$ has only one pre-image over each point of $[z,y[\subset S^\an(\CC)$,
\item[•]$f$ totally splits over $\pi_\CC^{-1}(]x,z[)$ (or over $\pi_\CC^{-1}([x,z[)$ if $\CC$ contains $x$).
\end{itemize}
\end{coro} 

\begin{proof}
If $I$ is an interval of $\R_{>0}$, denote by $\CC_I$ the analytic domain of $\mathbb{A}_k^{1,\an}$ defined by the condition $\vert T-1 \vert \in I$. This is an annulus. If $J$ is another interval of $\R_{>0}$, $\CC_I$ and $\CC_J$ are isomorphic if and only if $J\in \vert k^\times\vert\cdot I^{\pm 1}$. Therefore, as $z$ is a type-2 point, there exists some interval $I=]r_x, r_y[$ (or the same but closed in one or both extremities) of $\R_{>0}$, containing $p^{-\frac{p}{p-1}}$ in its interior, such that $\CC$ is isomorphic to $\CC_I$ through some isomorphism $\varphi : \CC\iso \CC_I$ that identifies $x, y$ and $z$ respectively with $\eta_{1,r_x}, \eta_{1, r_y}$ and $\eta_{1, p^{-\frac{p}{p-1}}}$.\\

Proposition \ref{décomposition du torseur sauvage} with $h=1$ tells us, for instance when $I=]r_x, r_y[$ is open, that the $\mu_p$-torsor $g : \G_m^\an \xrightarrow{z\mapsto z^{p}} \G_m^\an$ considered as a covering of $\CC_I$ ($f$ restricted to $g^{-1}(\CC_I)$) : 
\begin{itemize}
\item[•]totally splits over $\pi_{\CC_I}^{-1}\left( ]\eta_{1,r_x}, \eta_{1, p^{-\frac{p}{p-1}}}[\right)$,
\item[•]does not split with only one pre-image over $[\eta_{1, p^{-\frac{p}{p-1}}}, \eta_{1, r_y}[$, 
\end{itemize}
One conclude with $f=\varphi^* g\in H_1(\CC,\mu_p)$.
\end{proof}

\section{Construction of an admissible covering that is not overconvergent}

Let $\CC$ be a $k$-analytic annulus, defined as an analytic domain of $\A_{k}^{1,\an}$ by the condition $\lbrace \vert T\vert \in I\rbrace$, where $I$ is an interval of $\R_{>0}$ with non-empty interior. Let $x_0\in S^\an(\CC)$ that does not belong to the analytic boundary of $\CC$. There exists $r_0\in I$ such that $x_0=\eta_{0,r_0}$ (i.e. $x_0=r_0$ after identifying the skeleton $S^\an(\CC)$ with $I$). Consider the two subintervals $I^-, I^+\subset I$ defined by : $I^-=\lbrace i\in I, i\leqslant r_0\rbrace$ and $I^+=\lbrace i\in I, i \geqslant r_0\rbrace$. Let $\CC^-$ and $\CC^+$ the two sub-annuli of $\CC$ defined respectively by : $\CC^-=\lbrace \vert T\vert \in I^-\rbrace$ and $\CC^+=\lbrace \vert T\vert \in I^+\rbrace$. For any $r\in I$, let $\CC(r)$ be the circle of radius $r$, defined as a sub-annulus of $\CC$ by $\CC(r)=\lbrace \vert T\vert =r\rbrace$. The intersection of $\CC^-$ and $\CC^+$ corresponds to the circle of radius $r_0$ :
$\CC^-\cap \CC^+=\CC(r_0)$

\bigskip
As $k$ is non-trivially valued, from the density of type-$2$ points, there exist two sequences of type-$2$ points $(\alpha_n)_{n\in \N}\in S^\an(\CC^-)$ and $(\beta_n)_{n\in \N}\in S^\an(\CC^+)$ converging to $x_0$ and corresponding, after identification of $S^\an(\CC^\pm)$ with $I^\pm$, to an increasing sequence $(a_n)_{n\in \N}$ of $I^-$ and a decreasing sequence $(b_n)_{n\in \N}$ of $I^+$ both converging to $r_0$.

\bigskip
Let's apply Corollary \ref{corollary about splitting} to the annulus $\CC^-$ endowed with the type-$2$ point $\alpha_n$ : for each $n\in \N$, there exists a finite étale covering (even a $\mu_p$-torsor) $Y_n^-\to \CC^-$ that totally splits over $\pi^{-1}_{\CC^-}(]\alpha_n, x_0])$, and such that the fiber over each point of $S^\an(\CC^-)\setminus ]\alpha_n, x_0]$ has only one element. In particular, $Y_n^-$ totally splits over $\CC(r_0)=\pi^{-1}_{\CC^-}(\lbrace x_0\rbrace)$, and its restriction to $\CC(r_0-\varepsilon)$ is connected for each $\varepsilon>0$ such that $r_0-\varepsilon\in I^-$ and $n\gg 0$. Similarly, there exists a $\mu_p$-torsor $Y_n^+\to \CC^+$ that totally splits over $\pi^{-1}_{\CC^+}([x_0,\beta_n[)$, and such that the fiber over each point of $S^\an(\CC^+)\setminus [x_0, \beta_n[$ has only one element. In particular, $Y_n^-$ totally splits over $\CC(r_0)$, and its restriction to $\CC(r_0+\varepsilon)$ is connected for each $\varepsilon>0$ such that $r_0+\varepsilon\in I^+$ and $n\gg 0$.

\bigskip
We set $Y^-=\coprod_{n\in \N}Y_n^-$ and $Y^+=\coprod_{n\in \N}Y_n^+$, they are by construction objects of $\textbf{UFÉt}_{\CC^-}$ and $\textbf{UFÉt}_{\CC^+}$ respectively. Label the $p$ irreducible components of $Y^\pm\times_{\CC^\pm}\CC(r_0)$ by $\mathcal{Z}_{np+1}^\pm, \mathcal{Z}_{np+2}\pm, \ldots, \mathcal{Z}_{(n+1)p}^\pm$. Because of the splitting property of all our $\mu_p$-torsors over $\CC(r_0)$, each $\mathcal{Z}_{k}^\pm$ is mapped to $\CC(r_0)$ isomorphically and is a $k$-affinoid domain. Therefore, we can glue together $Y^-$ and $Y^+$ along their restrictions to $\CC(r_0)$, $Y^-\times_{\CC^-}\CC(r_0)=\coprod_{n\in N}\mathcal{Z}^-_n$ and $Y^+\times_{\CC^+}\CC(r_0)=\coprod_{n\in N}\mathcal{Z}^+_n$, by identifying $\mathcal{Z}^+_{n}$ with $\mathcal{Z}^-_{n+1}$, for all $n\in \N$.

\bigskip
This gives an étale morphism $Y\to \CC$ whose restrictions to $\CC^-$ and $\CC^+$ are respectively $Y^-\in \textbf{UFÉt}_{\CC^-}$ and $Y^+\in \textbf{UFÉt}_{\CC^+}$. As the affinoid domains $\CC^-$ and $\CC^+$ are open for the admissible topology (even though they are not open for the Berkovich overconvergent topology), the morphism $Y\to X$ is an object of $\textbf{Cov}_{\CC}^{\mathrm{adm}}$.

\begin{prop}
The morphism $Y\to \CC$ does not define an object of \emph{$\textbf{Cov}_{\CC}^{\mathrm{oc}}$}. Therefore, the inclusion \emph{$\textbf{Cov}_{\CC}^{\mathrm{oc}}\subseteq \textbf{Cov}_{\CC}^{\mathrm{adm}}$} is strict.
\end{prop}

\begin{proof}
If the morphism $Y\to \CC$ defined an object of $\textbf{Cov}_{\CC}^{\mathrm{oc}}$, there would exist an open overconvergent (i.e. for the Berkovich topology) neighbourhood $\mathcal{U}$ in $\CC$ of the Gauss point $\eta=\eta_{0,r_0}$ of $\CC(r_0)$ such that $Y\times_{\CC}\mathcal{U}\to \mathcal{U}$ would be a disjoint union of finite étale coverings of $\mathcal{U}$. In particular, the fiber over $\eta$ would be finite when restricted to any connected component of $Y\times_{\CC}\mathcal{U}$. It is then enough to show that there exists a connected component of $Y\times_{\CC}\mathcal{U}$ whose fiber over $\eta$ is infinite. \\
We can assume $\mathcal{U}^-=\mathcal{U}\cap \CC^-$ and $\mathcal{U}^+=\mathcal{U}\cap \CC^+$ to be connected. From \cite{ALY1}, Lemma $6.3.5$, there exists $\varepsilon>0$ such that $\mathcal{U}$ contains $\CC(r_0-\varepsilon)$ and $\CC(r_0+\varepsilon)$. Let $n_0\in \N$ be such that for all $n\geqslant n_0$, the restrictions of $Y_n^-$ and $Y_n^+$ to respectively $\CC(r_0-\varepsilon)$ and $\CC(r_0+\varepsilon)$ are connected. It implies that for all $n\geqslant n_0$, $Y^\pm_{n\;\vert\,\mathcal{U}}=Y_n^\pm\times_{\CC^\pm}\mathcal{U}^\pm$ is connected. Let $Y'$ be the image in $Y$ of the union of all the $Y_n^-$ and $Y_n^+$ for $n\geqslant n_0$. In the infinite sequence
$$(Y^-_{{n_0}\;\vert\,\mathcal{U}}, \;Y^+_{{n_0}\;\vert\,\mathcal{U}},\; Y^-_{{n_0+1}\;\vert\,\mathcal{U}}, \;Y^+_{{n_0+1}\;\vert\,\mathcal{U}},\; Y^-_{{n_0+2}\;\vert\,\mathcal{U}} \ldots ), $$

each term is a connected subset of $Y'_{\mathcal{U}}=Y'\times_{\CC}\,\mathcal{U}$, two consecutive ones have non-empty intersection, and the union of all of them is $Y'_{\mathcal{U}}$. This shows that $Y'_{\mathcal{U}}$ is connected, and the conclusion follows from the fact that $Y'_{\mathcal{U}}$ has an infinite fiber over $\eta$. 
\end{proof}

\begin{coro}
If $\overline{x}$ is a geometric point of $\CC$, the natural morphism of fundamental groups $$\pi_1^{\mathrm{dJ, \, adm}}(\CC, \overline{x})\to \pi_1^{\mathrm{dJ, \,oc}}(\CC,\overline{x}) $$ is not an isomorphism.
\end{coro}

\end{document}